\documentclass[11pt]{article}
\usepackage{amsfonts}
\usepackage{amssymb}
\usepackage{amssymb,amsfonts,amsmath,amsthm,cite, color}
\usepackage{epsfig}
\newtheorem{thm}{Theorem}[section]

\newtheorem{prop}{Proposition}[section]

\makeatletter \@addtoreset{equation}{section}

\def\mid{{\,|\,}}
\def\c3{{c\phi_3}}
\def\cc3{{\overline{c}\phi_3}}

\def\qed{\hfill \rule{4pt}{7pt}}

\def\Z{\mathbb Z}

\def\a{\alpha}

\begin{document}

\begin{center}
{{\Large\bf Congruences modulo powers of $5$ for  three-colored Frobenius partitions}}

\vskip 6mm

{ Xinhua Xiong
\\[%
2mm] Department of Mathmetics\\
Fudan University, Shanghai 200433,
P.R. China \\[3mm]
XinhuaXiong@fudan.edu.cn \\[0pt%
] }
\end{center}

\begin{abstract}
Motivated by a question of Lovejoy\,\,\cite{lovejoy}, we show that three-colored Frobenius partition
function $\c3$ and related arithmetic function $\cc3$ vanish modulo some powers of $5$ in certain arithmetic progressions.
To be more specific, we show that for every nonnegative integer $n$,
\begin{eqnarray*}
\begin{split}
\c3(45n+23)&\equiv 0\pmod{625},\\
\c3(45n+41)&\equiv 0\pmod{625},\\
\c3(75n+22)&\equiv 0\pmod{25},\\
\cc3(75n+72)&\equiv 0\pmod{25}.\\
\end{split}
\end{eqnarray*}
\end{abstract}



\section{Introduction}
The ordinary partition of a non-negative integer $n$ is a non-increasing sequence
of integers whose sum is $n$. A generalized Frobenius partition of $n$ is a two rowed array of
non-negative integers of the form
\begin{displaymath}
\left(
           \begin{array}{cccc}
             a_1 & a_2&\cdot\cdot\cdot & a_k \\
             b_1 & b_2&\cdot\cdot\cdot & b_k \\
           \end{array}
         \right),
\end{displaymath}
where the entries in each row are in non-increasing order and the integer $n$ that is partitioned
is $\sum_{i=1}^k(a_i+b_i+1).$ A $3$-colored Frobenius partition is an array of the above where the integer entries are
taken from $3$ distinct copies of the non-negative integers distinguished by color and the rows
ordered first by size and then by color with no two consecutive like entries in any row. \cite{Ono96} gave such an example.  If we denote by $p(n)$ the number of the ordinary partitions of $n$ and by $\c3(n)$ the $3$-colored Frobenius partitions of $n$ with the convention that $p(\a)=0$, if $\a \notin \Z$, Kotlitsch proved that the generating function for $ \cc3(n):=\c3(n)-p(n/3)$ is
\begin{displaymath}
\sum_{n = 0}^{\infty}\cc3(n)q^n = 9q \prod_{n \geq 1} \frac{(1-q^{9n})^3}{(1-q^{3n})(1-q^n)^3}.
\end{displaymath}
Since Ramanujan found his famous partition congruences,
their generalizations for various partition functions have been the subject of much investigation see\,\,\cite{Kolitsch89}\,\,\cite{Kolitsch90}\,\,\cite{Sellers}. Recently,
Ono\,\,\cite{Ono96} and Lovejoy\,\,\cite{lovejoy} proved the following congruences for small primes and every nonnegative integer $n$
\begin{eqnarray}
\begin{split}
\c3(45n+23)&\equiv 0\pmod{5},\\
\c3(45n+41)&\equiv 0\pmod{5},\\
\c3(63n+50)&\equiv 0\pmod{7},\\
\c3(99n+95)&\equiv 0\pmod{11},\\
\c3(171n+50)&\equiv 0\pmod{19}.\\
\end{split}
\end{eqnarray}

Moreover, Lovejoy asked the question: are there generalizations of the congruences (1.1) to
powers of $5$,\,$7$,\,$11$,\,$19$ analogous to the generalizations of Ramanujan's
congruences for the ordinary partition function? In this note, we prove

\begin{thm}\label{thm1.1} For every nonnegative integer $n$, we have
\begin{eqnarray*}
\begin{split}
\c3(45n+23)&\equiv 0\pmod{625},\\
\c3(45n+41)&\equiv 0\pmod{625}.\\
\end{split}
\end{eqnarray*}
\end{thm}
Also we give a theorem of another type of theorem 1.1:
\begin{thm}\label{thm1.2} For every nonnegative integer $n$, we have
\begin{eqnarray*}
\begin{split}
\c3(75n+22)&\equiv 0\pmod{25},\\
\cc3(75n+72)&\equiv 0\pmod{25}.\\
\end{split}
\end{eqnarray*}
\end{thm}
We note that there are no similar congruences for the arithmetic progressions and modulus listed in 1.1.

\section{Preliminaries}
When proving congruences we need to determine when the Fourier expansion of a modular form
has coefficients which are all multiples of $M$. So we define the $M$-adic order of a
formal power series.

\noindent Definition Let $M$ be a positive integer and $f = \sum_{n \geq N} a(n) q^n$ be a
formal power series in the variable $q$ with rational integer coefficients. The $M$-adic
order of $f$ is defined by
\begin{displaymath}
Ord_M(f) = \mbox{inf}\{n \mid a(n)\not\equiv 0\,\mbox{mod}\,M\}.
\end{displaymath}
Sturm proved the following criterion for determining whether two modular forms are
congruent for primes, Ono\,\,\cite{Ono96-2} noted the criterion holds for general integers as modulus.
\begin{prop}\label{prop2.1}
Let $M$ be a positive integer and  $f(z),\, g(z)\in M_k(\Gamma_0(N))$ with rational integers
satisfying
\begin{displaymath}
ord_M(f(z)-g(z))\ge 1+ \frac{kN}{12}\prod_{p}(1+\frac{1}{p}),
\end{displaymath}
 where the product is over
the prime divisors $p$ of $N$. Then $f(z)\equiv g(z)\ ({\rm mod\
}M)$, i.e., $ord_M(f(z)-g(z))=\infty$.
\end{prop}

Our proofs depend on special holomorphic modular forms on some congruence subgroups.
We use the following facts to construct such modular forms. Let
\begin{displaymath}
\eta(z):=q^{\frac{1}{24}}\prod_{n=1}^\infty (1-q^n)
\end{displaymath}
denote the Dedekind's Eta-function,
where $q=e^{2\pi iz}$ and ${\rm Im}(z)>0$. We know that it is a holomorphic modular
form of weight $\frac{1}{2}$ which does not vanish on complex upper half plane. A function
 $ f(z)$ is called an Eta-product if it can be written in the form of
\begin{displaymath}
f(z)=\prod_{\delta |N}\eta^{r_{\delta}}(\delta z),
\end{displaymath}
where $N$ and $\delta$ are positive integers and $r_{\delta}\in \mathbb{Z}$.

\begin{prop}\label{prop2.2} If $f(z)=\prod_{\delta|N}\eta^{r_{\delta}}(\delta
z)$ is an Eta-product satisfying
the following conditions:
\begin{equation*}\label{con1}
 k=\frac{1}{2}\sum_{\delta|N}r_{\delta}\in \mathbb{Z},
\end{equation*}
\begin{equation*}\label{con2}
\sum_{\delta|N}\delta r_{\delta}\equiv 0 \ ({\rm mod}\ 24),
\end{equation*}
\begin{equation*}\label{con3}
\sum_{\delta|N}\frac{N}{\delta} r_{\delta}\equiv 0 \ ({\rm mod}\
24),
\end{equation*}
then $f(z)$ satisfies
\begin{equation*}\label{relation1}
f\left(\frac{az+b}{cz+d}\right)=\chi(d)(cz+d)^kf(z)
\end{equation*}
for each $\left(
           \begin{array}{cc}
             a & b \\
             c & d \\
           \end{array}
         \right)\in \Gamma_0(N)$.
          Here the character $\chi$ is
         defined by $\chi(d):=\left(\frac{(-1)^ks}{d}\right)$, where
         \[ s:=\prod_{\delta|N}\delta^{r_{\delta}}\] and
         $\left(\frac{m}{n}\right)$ is Kronecker symbol.
\end{prop}
The analytic orders of an Eta-product at the cusps of $ \Gamma_0(N)$ was calculated by
Ligozat\,\,\cite{Ligozat75}.

\begin{prop}\label{prop2.3}
Let $c,d$ and $N$ be positive integers with $d|N$ and $(c,d)=1$. If
$f(z)$ is an Eta-product satisfying the conditions in Proposition
\ref{prop2.2} for $N$, then the order of vanishing of $f(z)$ at the
cusp $\frac{c}{d}$ is
\[
\frac{N}{24}\sum_{\delta
|N}\frac{(d,\delta)^2r_{\delta}}{(d,\frac{N}{d})d\delta}.
\]
\end{prop}

We also use the following proposition to construct modular forms, see Koblitz\,\,\cite{Koblitz84}.
\begin{prop}\label{prop2.4}
Let $\chi_1$ be a Dirichlet character modulo $N$,
and let $\chi_2$ be a primitive Dirichlet character modulo $M$.
Suppose $f(z)\in M_k(\Gamma_0(N),\chi_1)$ with Fourier expansion $$
f(z)=\sum\limits_{n=0}^\infty u(n)q^n.$$ Then
\noindent  for any positive integer $t|N$,\[f(z)|U(t):=\sum_{n=0}^\infty u(tn)q^n\]
  is the Fourier expansion of a modular form in $M_k(\Gamma_0(N),\chi_1)$.
\noindent  For any positive integer $t|N$,\[f(z)|V(t):=\sum_{n=0}^\infty u(n)q^{tn}\]
  is the Fourier expansion of a modular form in $M_k(\Gamma_0(tN),\chi_1)$.
\noindent  Let\[ g(z)=\sum\limits_{n=0}^\infty u(n)\chi_2(n)q^n. \] Then $g(z)\in
M_k(M^2N,\chi_1\chi_2^2)$.
\end{prop}

The last proposition we shall use is due to Eichhorn and Ono\,\,\cite{Eichhorn}.
\begin{prop}
If $l$ is a prime and $s\geq 1,$ then the Eta-product $\frac{\eta^{l^s}(z)}{\eta^{l^{s-1}}(lz)}$satisfies
\begin{eqnarray*}
\frac{\eta^{l^s}(z)}{\eta^{l^{s-1}}(lz)} &\equiv 1 \pmod {l^s}.
\end{eqnarray*}
\end{prop}

\section{The Congruence mod $625$}
In this section, we apply propositions above to prove Theorem 1.1.
\begin{thm}\label{thm3.1}
Define $a(n)$ by the infinite product
\begin{displaymath}
\sum_{n=0}^{\infty} a(n) q^n = \prod_{n=1}^{\infty} \frac {1}{(1-q^{3n})(1-q^n)^3},
\end{displaymath}
then the coefficients of $a(n)$ satisfies
\begin{displaymath}
a(n) \equiv 0 \pmod {625}\,\mbox{if}\,\, n \equiv 13,\,22,\,31,\,40 \pmod{45}.
\end{displaymath}
\end{thm}

\noindent {\it Proof.}\,Define the following Eta-product
\begin{eqnarray*}
g(z)&=&\frac{\eta^{13}(45z)\eta^3(135z)}{\eta(3z)\eta^3(z)}
\left(\frac{\eta^{625}(z)}{\eta^{125}(5z)}\right)^2\nonumber\\[5pt]
&=& q^{41}\left(\prod_{n=1}^{\infty}\frac{1}{(1-q^{3n})(1-q^n)^3}\right)\\[5pt]
& &\cdot \left(\prod_{n=1}^{\infty}(1-q^{45n})^{13}(1-q^{135n})^3\right)\\[5pt]
& &\cdot\left(\frac{\prod_{n=1}^{\infty}(1-q^n)^{1250}}
{\prod_{n=1}^{\infty}(1-q^{5n})^{250}}\right)\nonumber\\[5pt]
&:=& \sum_{n\geq 41} r(n)q^n.
\end{eqnarray*}
By Propositions 2.2 and 2.3 it turns out that $g(z)\in S_{506}(\Gamma_0(135),Id)$,
where $Id$ is the trivial Dirichlet character mod $135$. By Proposition 2.5 the product of the first $3$ factors of the middle expression
above possesses the congruence properties of the Fourier coefficients of $g(z)$ mod $625$.
We note that  Theorem 3.1 is equivalent to the congruences
\begin{displaymath}
r(n)\equiv 0 \pmod{625} \,\mbox{if}\, n\equiv 9,\,18,\,27,\,36 \pmod{45}.
\end{displaymath}
Use the Proposition 2.4 we find that the modular form
\begin{eqnarray*}
& &g(z)\mid U(9) - g(z)\mid U(45)V(5)\\[5pt]
&=& \sum r(9n)q^n -\sum r(45n)q^{5n}\nonumber\\[5pt]
&=& \sum_{n\not\equiv 0\pmod{5}} r(9n)q^n
\end{eqnarray*}
is in $M_{506}(\Gamma_0(675),Id).$ By Sturm's criterion, if it can be shown that $r(9n)\equiv 0
\pmod{625}$ when $ n\leq 45541$ and $n\not\equiv 0 \pmod{5}$, then our theorem follows. These have been
verified by machine computation so we proved the theorem.

\noindent {\it Proof of Theorem \ref{thm1.1}.} By Jacobi's triple product identity, we have
\begin{eqnarray*}
&&\c3(45n+23)\\
&=& \cc3(45n+23)\\
           &= &9\sum_{k\geq 0}(-1)^k(2k+1)(a(45n+23-(1+\frac{9k^2+9k}{2}))),\\
&&\c3(45n+41)\\
&=& \cc3(45n+41)\\
           &= &9\sum_{k\geq 0}(-1)^k(2k+1)(a(45n+41-(1+\frac{9k^2+9k}{2}))).
\end{eqnarray*}
Since modulo $45$, $1+\frac{9k^2+9k}{2}$ is $1,\,10,\,28$, so $45n+23-(1+\frac{9k^2+9k}{2})$ and
$45n+41-(1+\frac{9k^2+9k}{2})$ are congruence to $13,\,22,\,31,\,40$ modulo $45$. By Theorem 3.1, Theorem 1.1 is proved.

\section{The Congruence mod $25$}
In this section we prove Theorem 1.2.

\noindent {\it Proof of Theorem \ref{thm1.2}.}  Define the following Eta-product
\begin{eqnarray*}
f(z)&=&9\frac{\eta^{3}(9z)\eta(75z)}{\eta(3z)\eta^3(z)}
\left(\frac{\eta^{25}(z)}{\eta^{5}(5z)}\right)^2\nonumber\\[5pt]
&=& 9q^{4}\left(\prod_{n=1}^{\infty}\frac{(1-q^{9n})^3}{{(1-q^{3n})(1-q^n)^3}}\right)
\left(\prod_{n=1}^{\infty}(1-q^{75n})\right)\\
&&\cdot\left(\prod_{n=1}^{\infty}\frac{(1-q^n)^{50}}{(1-q^{5n})^{10}}\right)\nonumber\\[5pt]
&=& q^3\left(\sum_{n=0}^{\infty}\cc3(n)q^n\right)\left(\prod_{n=1}^{\infty}(1-q^{75n})\right)\left(\prod_{n=1}^{\infty}\frac{(1-q^n)^{50}}
{(1-q^{5n})^{10}}\right)\nonumber\\[5pt]
&:=& \sum_{n\geq 4} c(n)q^n.
\end{eqnarray*}
By Propositions 2.2 and 2.3, We find $f(z)$ is in $S_{20}(\Gamma_0(225),Id)$, where $Id$ is the trivial character
modulo $225.$  We note that our Theorem 1.2 is equivalent to the following congruences: For every nonnegative integer $n$
\begin{eqnarray}
\begin{split}
c(75n+25)&\equiv 0\pmod{25},\\
c(75n)&\equiv 0\pmod{25}.\\
\end{split}
\end{eqnarray}
Let
\begin{displaymath}
f_2(z)=\sum_{n\geq 4}c(n)\left(\frac{n}{3} \right) q^n.
\end{displaymath}
By Proposition 2.4, $f_2(z)$ is in $S_{20}(\Gamma_0(2025),Id)$. Define another modular form
$F(z)\in S_{20}(\Gamma_0(2025),Id)$ by
\begin{eqnarray*}
F(z)&:=& \sum_{n=0}^{\infty}d(n)q^n\\
 &= &f(z)+f_2(z)\\
&=&\sum_{\left(\frac{n}{3} \right)=1}2c(n)q^n +\sum_{n\equiv 0\,\mbox{mod}\,3}c(n)q^n.
\end{eqnarray*}
Apply Hecke operator $U(25)$, we obtain
\begin{displaymath}
F(z)\mid U(25)= \sum_{n\geq 1}d(75n)q^n.
\end{displaymath}
If $F(z)\mid U(25)\equiv 0 \pmod{25},$ then 4.1 holds. By
Sturm's criterion, we need to check that the congruence holds for the first $25\cdot(20/12)[SL_2(Z):\Gamma_0(2025)]+1
=135001$ terms which is easily verified.
\qed

\section{Conclusion}
We find there are no similar congruences for powers of \,$7,\,\,11,\,\,19$ \,in arithmetic
 progressions given in 1.1. So we ask: are there congruences for higher powers of five or other primes in the arithmetic progressions of type $3^\a5^{\beta}\cdot n+\gamma(\a,\beta)$? Second, find identities for $\c3$, which interpret congruences above for $\c3$, analogous to the Ramanujan-type identities for the ordinary parttition function interpret Ramanujan's congruences. We note that similar results were obtained by Sellers\,\,\cite{Sellers}.

\vspace{.2cm} \noindent{\bf Acknowledgments.} The author wish to thank Professor Meng Chen for helpful comments. The author also would like to thank the referee for his/her helpful corrections and suggestions.


\end{document}